\documentclass{amsart}  %
\usepackage{latexsym,amssymb,amsmath,amsfonts,amsthm}
\usepackage{tabularx}
\newcommand{\bC}{{\mathbb C}}
\newcommand{\bR}{{\mathbb R}}
\newcommand{\bN}{{\mathbb N}}
\newcommand{\sF}{{\mathcal{F}}}
\newcommand{\sS}{{\mathcal{S}}}

\newcommand{\norm}[1]{\left\Vert #1\right\Vert}
\newcommand{\bignorm}[1]{\bigl\Vert #1\bigr\Vert}
\newcommand{\Bignorm}[1]{\Bigl\Vert #1\Bigr\Vert}

\newcommand{\abs}[1]{\left\vert #1\right\vert}
\newcommand{\bigabs}[1]{\bigl\vert #1\bigr\vert}
\newcommand{\Bigabs}[1]{\Bigl\vert #1\Bigr\vert}
\newcommand{\biggabs}[1]{\biggl\vert #1\biggr\vert}

\newcommand{\Mu}{\mathrm{M}}
\def\means{\mathbin{``{=}"}}

\newtheorem{theorem}{Theorem}

\newtheorem*{theorem*}{Theorem}

\theoremstyle{remark}
\newtheorem{remark}[theorem]{Remark}

\begin{document}
\title[Hardy's uncertainty principle for operators]{Vector-valued distributions and \\
Hardy's uncertainty principle for operators}

\author{M. G.~Cowling}
\address{School of Mathematics\\ University of Birmingham\\ Edgbaston Birmingham B15 2TT\\ UK}
\email{M.G.Cowling@bham.ac.uk}

\author{B.~Demange}
\address{UFR de Math\'ematiques\\ Universit\'e de Grenoble\\ BP 74\\ 38402 St Martin d'H\`eres, France}
\email{Bruno.Demange@ujf-grenoble.fr}

\author{M.~Sundari}
\address{Chennai Mathematical Institute\\
Plot H1 SIPCOT IT Park\\
Padur P O\\
SIRUSERI 603 103}
\email{sundari@cmi.ac.in}

\date{May 6, 2008}

\maketitle

\begin{abstract}
Suppose that $f$ is a function on $\mathbb{R}^n$ such that
$\exp(a \abs{\,\cdot\,}^2) f$ and $\exp(b \abs{\,\cdot\,}^2)\hat f$ are bounded, where $a,b>0$.  Hardy's Uncertainty Principle asserts that if $ab > \pi^2$, then $f = 0$,
while if $ab = \pi^2$, then $f = c \exp(-a\abs{\,\cdot\,}^2)$.  In this paper,
we generalise this uncertainty principle to vector-valued functions, and hence to operators.  The principle for operators can be formulated loosely by saying that the kernel of an operator cannot be localised near the diagonal if the spectrum is also localised.
\end{abstract}

\section{Introduction}

We define the Fourier transform $\hat f$ of a function $f$ on $\bR^n$ by
\begin{equation}\label{def-Fourier}
\hat f (\xi) := \int_{\bR^n} f(x) \exp(-2  \pi i \xi \cdot x) \,dx
\qquad\forall \xi \in \bR^n,
\end{equation}
provided that this makes sense.

G.H. Hardy \cite{Hardy33} showed that if $f$ is a function on $\bR^n$ such that
$$\begin{aligned}
\bigabs{ f(x) } &\leq c \exp(-a \abs{x}^2) \qquad\forall x \in \bR^n \\
\bigabs{ \hat f(\xi) } &\leq c \exp(-b \abs{\xi}^2) \qquad\forall \xi \in \bR^n ,
\end{aligned}$$
where $a,b > 0$, and if
$ab > \pi^2$, then $f = 0$, while if  $ab = \pi^2$, then $f(x) = c' \exp(-a \abs{x}^2)$.  Clearly the second result implies the first.

By rescaling, we may and shall suppose without loss of generality that $a=\pi$ in  the above results.  Similarly we may remove the constant $c$.

For a detailed account of the proof and some classical related results, see H. Dym and H.P. McKean \cite{DyMc72}.

This result, which we call Hardy's Uncertainty Principle, has been generalised in many ways.  For instance, G.W. Morgan \cite{Morgan34} obtained a sharp version involving the modified exponentials $\exp(-a \abs{\,\cdot\,}^p)$ and $\exp(-b\abs{\,\cdot\,}^q)$,
where $1/p + 1/q  =1$, giving the extremal functions.
Cowling and J.F. Price \cite{CoPr84}
showed that if $p$ or $q$ is finite, and
\[\begin{aligned}
\exp(\pi \abs{\,\cdot\,}^2)f &\in L^p(\bR^n)  \\
\exp(\pi \abs{\,\cdot\,}^2) \hat f &\in L^q(\bR^n)  ,
\end{aligned}\]
then $f =0$.

Demange (see A. Bonami and Demange \cite{BoDe06}) proved the best version for the usual quadratic exponentials.

\begin{theorem}\label{Demangethm}
Suppose that $\Phi \in \sS(\bR^n)'$ (that is, $\Phi$ is a tempered distributions) and
\[\begin{aligned}
\exp(\pi \abs{\,\cdot\,}^2)\Phi &\in \sS(\bR^n)'  \\
\exp(\pi \abs{\,\cdot\,}^2) \hat \Phi &\in \sS(\bR^n)'  ,
\end{aligned}\]
then
\[
\Phi(x) = \sum_{\abs{\alpha} \leq N} c_\alpha \, x^\alpha \, \exp(-\pi \abs{x}^2)
\qquad\forall x \in \bR^n,
\]
that is, $\Phi$ is the distribution arising by integration against the function on the right hand side.
\end{theorem}

There has been interest recently in developing our understanding of heat diffusion and its relationship with uncertainty principles into more general contexts than $\bR^n$, such as the Heisenberg group, where the work of Thangavelu \cite{Thang} shows that the subelliptic world presents some surprises, Lie groups \cite{VSC}, where structure plays an important role, and differential operators on measured metric spaces \cite{CRS}.

In this paper, we extend the uncertainty principles of Hardy and of Demange
to vector-valued functions.  The former generalises naturally to Banach space valued functions, while the latter can be formulated for functions with values in the dual of a
Fr\'echet space.

We then consider uncertainty principles for operators.  Our main theorems are concerned with operators from $\sS(\bR^n)$ to $\sS(\bR^n)'$, where the kernel of the operator and one of its partial Fourier transforms satisfy conditions like those of Theorem~\ref{Demangethm}, and with operators on $L^2(\bR^n)$, whose kernels are locally integrable and which satisfy inequalities like those of Hardy's theorem.
We extend and sharpen previous work of Cowling and Sundari~\cite{CS} which led us to conjecture the result proved here.  
All the results for operators are stated and proved in Section~\ref{s:Operators}.

Here is one of our main results, in which $P_t$ is the heat operator and $p_t$ is its kernel.

\begin{theorem}\label{main-theorem-integrals}
Suppose that $t$ is in $\bR^+$, that $K$ is the operator on $L^2(\bR^n)$ associated to a locally integrable kernel $k$, and that
$$\begin{aligned}
\bigabs{ k(x,y) }   &\leq  p_t(x,y)
\qquad\forall x,y \in \bR^n \\
\bignorm{Kf}_{2}  & \leq  \bignorm{P_t f}_{2}
\qquad\forall f \in L^2(\bR^n).
\end{aligned}$$
Then there is a bounded measurable function $m$ on $\bR^n$ such that
\[
k(x,y) = m(x) \, p_t(x,y)
\qquad\forall x,y \in \bR^n,
\]
and the operator $K$ is the the heat operator $P_t$ followed by multiplication by $m$.
\end{theorem}

Our other main result, Theorem \ref{thm:BDJ-op}, generalises Theorem \ref{main-theorem-integrals} in the same way that Theorem \ref{Demangethm}  generalises Hardy's original result.

In the next section, we summarise L. Schwartz' theory of distributions, including the case of vector-valued distributions, and establish versions of Hardy's Uncertainty Principle in this context.  In Section \ref{s:Operators}, we consider the theory for linear operators and prove Theorems 1 and 4.

\section{Tempered distributions}
\subsection{The Schwartz space}
Recall that the Schwartz space $\sS(\bR^n)$ is defined to be the space of all smooth functions $f\colon\bR^n\to\bC$ such that $\norm{f}_{(N)} < \infty$ for all $N$ in $\bN$, where
\begin{equation}\label{e:S-norms}
\norm{f}_{(N)} := \biggl(\sum_{\abs{\alpha}, \abs{\beta} \leq N} \bignorm{ {}\cdot{}^\alpha D^\beta f}_{2}^2\biggr)^{1/2} .
\end{equation}
Here we use standard multi-index notation.  Thus, when $\alpha \in \bN^n$, we define
$x^\alpha$ to be $x_{1}^{\alpha_{1}} \dots x_n^{\alpha_n}$, and
\[
D^\beta f(x) := \frac{\partial^{\beta_{1}}}{\partial x_{1}^{\beta_{1}}} \cdots \frac{\partial^{\beta_n}}{\partial x_n^{\beta_n}} f(x) .
\]

Then $\sS(\bR^n)$ is the intersection of the completions $\sS_{(N)}(\bR^n)$ of the spaces of compactly supported smooth functions in the norm $\norm{\cdot}_{(N)}$.  It is a
Fr\'echet space: it has a metric,  for instance,
\[
d(f,g) := \sum_{N \in \bN} \,\frac{\norm{f-g}_{(N)}}{2^{N} (1+ \norm{f-g}_{(N)})}
\qquad\forall f,g \in \sS(\bR^n),
\]
but in general $d(\lambda f, \lambda g) \neq \abs{ \lambda }  d(f,g)$.

\subsection{Tempered distributions}
The space of tempered distributions  $\sS(\bR^n)'$ is the dual of $\sS(\bR^n)$.
It is the union of the spaces  $\sS_{(N)}(\bR^n)'$.
Thus, any tempered distribution $\Phi$ belongs to some space $\sS_{(N)}(\bR^n)'$, and then
\[
|\Phi(f)| \leq c \norm{f}_{(N)} \qquad\forall f \in \sS(\bR^n).
\]

Every slowly growing locally integrable function $k$, that is, a locally integrable functions whose total mass over a ball grows polynomially with the radius of the ball, defines a distribution $\Phi_k$ by integration:
\[
\Phi_k(f) := \int_{\bR^n} k(x) \,f(x) \,dx.
\]
A statement of the form $\Phi = f$, where $\Phi$ is a distribution and $f$ is a function, is to be interpreted that $\Phi$ is the distribution obtained by integrating against $f$.
%
%
%
For many purposes, we may consider all distributions as being given by integration against a function, as they may be derived from such distributions by standard analytical processes, such as taking limits or differentiating.

Sometimes we write a formula involving distributions using pointwise notation, as shorthand for the integrated version of the formula.  Thus, in $\bR^{2n}$, we will write
\[
\Phi(x,y) = \Mu(x) \, g(x-y)
\qquad\forall x,y \in \bR^n,
\]
where $g\in \sS(\bR^n)$ and $\Mu \in \sS(\bR^n)'$.
This ``means'' that
\[\begin{aligned}
\Phi(f)
& \means \int_{\bR^{n}} \int_{\bR^{n}} \Mu(x) \, g(x-y)  \,f(x,y) \,dy \,dx \\
& \means \int_{\bR^{n}} \Mu(x) \int_{\bR^{n}} g(x-y)  \,f(x,y) \,dy \,dx
\end{aligned}\]
for all $f$ in $\sS(\bR^{2n})$, and so really means that applying the distribution $\Phi$ to $f$ is the same as applying the distribution $\Mu$ to the function
$x \mapsto \int_{\bR^n} g(x-y) \, f(x,y) \,dy$.

Every smooth function $h$, all of whose partial derivatives are of polynomial growth, multiplies $\sS(\bR^n)$ pointwise, and hence multiplies $\sS(\bR^n)'$ by duality: for $\Phi$ in $\sS(\bR^n)'$,
\[
h\Phi(f) := \Phi(hf)
\qquad\forall f \in \sS(\bR^n)  .
\]
In particular, we write $E_{\rm i,1}$ for the function $(x,y) \mapsto \exp(2\pi i x\cdot y)$ on $\bR^{2n}$.  Then $E_{\rm i,1}$ multiplies $\sS(\bR^{2n})$ and hence $\sS(\bR^{2n})'$.

For $f, g$ in $\sS(\bR^n)$, we define $f \otimes g$ in $\sS(\bR^{2n})$ by
\[
f \otimes g(x,y) := f(x) \, g(y)
\qquad\forall x,y \in \bR^n.
\]
The set of finite linear combinations of such ``outer tensor product functions'' is dense in $\sS(\bR^{2n})$.

The Schwartz kernel theorem states that there is a one-to-one correspondence between continuous linear operators $T$ from $\sS(\bR^n)$ to $\sS(\bR^n)'$ and tempered distributions $\Phi$ in $\sS(\bR^{2n})'$, described by the formula
\[
T(f)(g) = \Phi(g \otimes f)
\qquad\forall f, g \in \sS(\bR^n).
\]
At least formally, we may write $T(f)(x) = \int_{\bR^n} \Phi(x,y) \, f(y) \,dy$.

\subsection{The Fourier transformation}
The Fourier transformation $\sF\colon f \mapsto \hat f$ (where $\hat f$ is defined by \eqref{def-Fourier}) is a bijection of $\sS(\bR^n)$.  Generalising the formula
\[
\int_{\bR^n} f(x) \, \hat g (x) \,dx =  \int_{\bR^n} \hat f(\xi) \, g (\xi) \,d\xi
\qquad\forall f,g \in L^1(\bR^n),
\]
we define the Fourier transform $\hat\Phi$ of a tempered distribution $\Phi$:
\[
\hat\Phi(f) := \Phi(\hat f)
\qquad\forall f\in \sS(\bR^n) .
\]
This extended Fourier transformation is also bijective on $\sS(\bR^n)'$.

Suppose that $f$ is a function on $\bR^{2n}$.  The partial Fourier transform $\sF_{2} f$ of $f$ is defined by
\[\begin{aligned}
\sF_{2} f (x,\eta) &:= \int_{\bR^n} f(x,y) \, \exp(-2  \pi i \eta \cdot y) \,dy
\qquad\forall x, \eta \in \bR^n
\end{aligned}\]
provided that this makes sense.  It is easy to show that the partial Fourier transformation $\sF_{2}$ is a bijection of $\sS(\bR^{2n})$, and hence extends to $\sS(\bR^{2n})'$.

\subsection{Hardy's Uncertainty Principle for Distributions}

As already mentioned, Demange (see \cite{BoDe06}) proved the best version of Hardy's result for the usual quadratic exponentials:  if $\Phi \in \sS(\bR^n)'$  and
\[\begin{aligned}
\exp(\pi \abs{\,\cdot\,}^2)\Phi &\in \sS_{(N)}(\bR^n)'  \\
\exp(\pi \abs{\,\cdot\,}^2) \hat \Phi &\in \sS_{(N)}(\bR^n)'  ,
\end{aligned}\]
then
\[
\Phi(x) = \sum_{\abs{\alpha} \leq N'} c_\alpha \, x^\alpha \, \exp(-\pi \abs{x}^2)
\qquad\forall x \in \bR^n.
\]

We claim that we may assume that $N'\leq N$.  To see this, we may suppose that $c_\alpha \neq 0$ for some $\alpha$ such that $\abs{\alpha} = N'$, in which case, for some $y_0$ in $\bR^n$, there exists a positive constant $c$ such that
$$
\bigabs{\sum_{\abs{\alpha}=N'} c_\alpha (ty_0)^\alpha} = c t^{N'}
\qquad\forall t \in \bR^+.
$$
Now if $\Phi$ satisfies the hypotheses and hence also the conclusion of the theorem, then
\[\begin{aligned}
\biggabs{ \sum_{\abs{\alpha} \leq N'} c_\alpha \int_{\bR^n} {x}^\alpha \, f(x) \,dx }
&= \bigabs{\Phi\bigl( \exp(\pi \abs{\cdot}^2) f \bigr)}
= \bigabs{\bigl(\exp(\pi \abs{\cdot}^2)\Phi\bigr) (  f )} \\
&\leq c \norm{f}_{(N)}
\qquad\forall f \in \sS(\bR^n) .
\end{aligned}\]
Take a smooth function $f$ with small compact support such that $\int_{\bR^n} f(x) \,dx =1$.  Then, on the one hand,
\[
\norm{\tau_y f}_{(N)} \leq c\, (1 + \abs{y})^{N}
\qquad\forall y \in \bR^n ,
\]
where $\tau_y f$ denotes the translate $f({\cdot}-y)$ of $f$, while on the other,
\[\begin{aligned}
\Bigabs{\sum_{\abs{\alpha} \leq N'} c_\alpha \int_{\bR^n} {x}^\alpha \, \tau_{y} f(x) \,dx }
&=\Bigabs{\sum_{\abs{\alpha} \leq N'} c_\alpha \int_{\bR^n} (x+y)^\alpha \, f(x) \,dx } \\
&= \Bigabs{\sum_{\abs{\alpha} = N'} c_\alpha \, {y}^\alpha }
+ O(|y|^{N'-1}),
\end{aligned}\]
and the claim follows by taking $y$ to be $ty_0$ and letting $t$ grow.

\subsection{Vector-valued distributions}

Take a Banach space $X$.

For nice enough $X$-valued functions $f$, we may still compute partial derivatives and multiply by scalar-valued functions.  Hence we may form the vector-valued Schwartz space $\sS(\bR^n;X)$ of $X$-valued functions $f$ such that $\norm{f}_{(N)}<\infty$, for all $N$ in $\bN$, where the norm $\norm{\cdot}_{(N)}$ is still defined by \eqref{e:S-norms}.

It is then possible to define the Fourier transform $\hat f$ of such a function, by formula~ \eqref{def-Fourier}, where the integrand is $X$-valued, and to prove an $X$-valued version of Hardy's Uncertainty Principle.

\begin{theorem}\label{vector-valued-thm1}
Suppose that $f\colon\bR^n\to X$ is a smooth vector-valued function such that
\[\begin{aligned}
\bignorm{f(x)}_{X} &\leq \exp(-\pi \abs{x}^2) \qquad\forall x \in \bR^n \\
\bignorm{\hat f(\xi)}_{X} &\leq \exp(-\pi \abs{\xi}^2) \qquad\forall \xi \in \bR^n .
\end{aligned}\]
Then $f(x) = f(0) \exp(-\pi \abs{x}^2)$.
\end{theorem}

\begin{proof}
Take an element $V$ of the dual space $X'$, and denote by $V(f)$ the scalar-valued function $x \mapsto V(f(x))$.  Clearly this satisfies the inequalities
\[\begin{aligned}
\bigabs{V(f)(x)} &\leq c \exp(-\pi \abs{x}^2) \qquad\forall x \in \bR^n \\
\bigabs{(V(f))\hat {\phantom{f}}(\xi)} &\leq c \exp(-\pi \abs{\xi}^2) \qquad\forall \xi \in \bR^n ,
\end{aligned}\]
so $V(f)(x) = c' \exp(-\pi \abs{x}^2)$ for all $x$ in $\bR^n$.

Evidently $c'$ depends linearly on $V$, and so there is an element $C$ of $X''$ such that
$c' = C(V)$.  Now, by taking $x$ to be $0$, we see that $C(V) = V(f(0))$.
\end{proof}

Take a Fr\'echet space $X$, which is the intersection of a decreasing family of Banach spaces $X_{(M)}$ with increasing norms $\norm{\cdot}_{(M)}$.   The standard
Fr\'echet metric $d$ is defined by
\[
d(u,v) := \sum_{M \in \bN} \,\frac{\norm{u-v}_{(M)}}{2^{n} (1+ \norm{u-v}_{(M)})}
\qquad\forall u,v \in X.
\]

Suppose that $f$ is $X$-valued and continuous, and that for some $M$ in $\bN$,
$$\begin{aligned}
\bignorm{f(y)}_{(M)} &\leq \exp(-\pi \abs{y}^2) \qquad\forall y \in \bR^n \\
\bignorm{\hat f(\eta)}_{(M)} &\leq \exp(-\pi |\eta|^2) \qquad\forall \eta \in \bR^n .
\end{aligned}$$
Then $f(y) = f(0) \exp(-\pi \abs{y}^2)$, by Theorem \ref{vector-valued-thm1}.

If $d$ is the Fr\'echet metric just defined, and $f\colon \bR^n \to X$ is continuous and
\[\begin{aligned}
d(f(y),0) &\leq c \exp(-\pi \abs{y}^2) \qquad\forall y \in \bR^n \\
d(\hat f(\eta),0) &\leq c \exp(-\pi |\eta|^2) \qquad\forall \eta \in \bR^n ,
\end{aligned}\]
then the same conclusion holds.  Indeed, assuming (as we may) that $c \geq 1$, then these estimates imply that
\[\begin{aligned}
\frac{\bignorm{f(y)}_{(0)}}{1+\norm{f(y)}_{(0)}} &\leq c \exp(-\pi \abs{y}^2) \qquad\forall y \in \bR^n \\
\frac{\bignorm{\hat f(\eta)}_{(0)}}{1+\bignorm{\hat f(\eta)}_{(0)}} &\leq c \exp(-\pi |\eta|^2) \qquad\forall \eta \in \bR^n ,
\end{aligned}\]
 and hence that, when $y$ and $\eta$ are big enough so that $c \exp(-\pi \abs{y}^2) < 1/2$ and $c \exp(-\pi |\eta|^2) < 1/2$,
\[\begin{aligned}
\bignorm{f(y)}_{(0)} &\leq 2c \exp(-\pi \abs{y}^2) \\
\bignorm{\hat f(\eta)}_{(0)} &\leq 2c \exp(-\pi |\eta|^2)  .
\end{aligned}\]
It follows that $y \mapsto V'(f(y))$ is a continuous function on $\bR^n$ for any $V'$ in ${X_{(0)}}'$, and that
\[\begin{aligned}
\abs{V'(f(y))} &\leq 2c \norm{V'}_{(0)} \exp(-\pi \abs{y}^2)
\qquad\forall y \in \bR\setminus [-\kappa,\kappa]\\
\abs{V'(\hat f(\eta))} &\leq 2c \norm{V'}_{(0)} \exp(-\pi |\eta|^2)
\qquad\forall \eta \in \bR\setminus [-\kappa,\kappa] ,
\end{aligned}\]
where $\kappa = \bigl( \log(2c) / \pi \bigr)^{1/2}$.  Thus,  $\exp(\pi \abs{\,\cdot\,}^2) V'(f(\cdot))$ and $\exp(\pi \abs{\,\cdot\,}^2) V'(\hat f(\cdot))$ are distributions, and by Theorem \ref{Demangethm},
\[
V'(f(y)) = c(V') \exp(-\pi \abs{y}^2)
\qquad \forall y \in \bR^n,
\]
where $c(V')$ is a constant depending on $V'$.  But $c(V') = V'(f(0))$ for all $V'$ and hence
\[
f(y) = f(0) \exp(-\pi \abs{y}^2)
\qquad \forall y \in \bR^n.
\]
Note that there are metrics (e.g., $d^{1/2}$) so that Hardy's Uncertainty Principle cannot be formulated as we have formulated the result for $d$.  Thus this result is more a curiosity than ``a theorem from the book''.\\[4pt]

We may generalise further, and consider distributions with values in $X'$, the dual of the Fr\'echet space $X$.  At least formally, every element $\Phi$ of $\sS(\bR^n ;X)'$ is such an object: for any such $\Phi$ there exist $M$ and $N$ in $\bN$ such that
\[
|\Phi(f)| \leq c \norm{f}_{(M,N)}
\qquad\forall f \in \sS(\bR^n;X),
\]
where
\[
\norm{f}_{(M,N)} := \biggl(\sum_{\abs{\alpha}, \abs{\beta} \leq N} \Bignorm{ \bignorm{{}\cdot{}^\alpha D^\beta f }_{(M)} }_{2}^2\biggr)^{1/2} .
\]

Given an $X'$-valued distribution $\Phi$ and an element $V$ of $X$, we may define a scalar-valued distribution $\Phi(V)$ by the formula
\[
\Phi(V)(f) := \Phi (fV)
\qquad\forall f \in \sS(\bR^n),
\]
where $fV$ is the function whose value at $x$ is $f(x)V$.
Thus defined, $\Phi(V)$ is a distribution, because for $V$ in $X$ and a scalar-valued function $f$ in $\sS(\bR^n)$, the vector-valued function $fV$ is in $\sS(\bR^n;X)$.

\subsection{A vector-valued version of Theorem \ref{Demangethm}}

\begin{theorem}\label{thm:BDJ-vv}
Suppose that $t$ is in $\bR^+$, that $X$ is a Fr\'echet space, that $\Phi$ is an $X'$-valued distribution, and that
\[\begin{aligned}
\exp(t \pi \abs{\,\cdot\,}^2)\Phi &\in \sS_{}(\bR^n;X)'  \\
\exp(\frac{\pi}{t} \abs{\,\cdot\,}^2) \hat \Phi &\in \sS_{}(\bR^n;X)'  .
\end{aligned}\]
Then
\[
\Phi(y) = \sum_{\abs{\alpha} \leq N} C_\alpha \, y^\alpha \, \exp(-t \pi \abs{y}^2)
\qquad\forall y \in \bR^n,
\]
where $C_\alpha \in X'$.
\end{theorem}

In particular, this applies when $X = \sS(\bR^n)$.

\begin{proof}
By rescaling if necessary, we may suppose without loss of generality that $t=1$.

As already remarked, there exist $M$ and $N$ in $\bN$ such that
\[\begin{aligned}
\exp(\pi \abs{\,\cdot\,}^2)\Phi &\in \sS_{(M,N)}(\bR^n;X)'  \\
\exp(\pi \abs{\,\cdot\,}^2) \hat \Phi &\in \sS_{(M,N)}(\bR^n;X)'  .
\end{aligned}\]
For $V$ in $X$, it follows that
\[\begin{aligned}
\exp(\pi \abs{\,\cdot\,}^2)\Phi(V) &\in \sS_{(N)}(\bR^n)'  \\
\exp(\pi \abs{\,\cdot\,}^2) \hat \Phi(V) &\in \sS_{(N)}(\bR^n)'  ,
\end{aligned}\]
and by Theorem \ref{Demangethm},
\[
\Phi(V)(y) = \sum_{\abs{\alpha} \leq N}  c_\alpha(V) \, y^\alpha \, \exp(-\pi \abs{y}^2)
\qquad\forall y \in \bR^n.
\]
We now determine how the numbers $c_\alpha(V)$ depend on $V$.

The functions $y \mapsto {y}^\alpha \,\exp(-\pi \abs{y}^2)$ are linearly independent for different $\alpha$, and we may find Schwartz functions $f_\beta$ such that
\[
\int_{\bR^n} {y}^\alpha \,\exp(-\pi \abs{y}^2)\,f_\beta(y) \,dy  = \delta_{\alpha,\beta}
\]
(the Kronecker delta) for all $\alpha$ and $\beta$ such that $\abs{\alpha}, \abs{\beta} \leq N$.  Now $c_\beta(V)$ is equal to $\Phi(V)(f_\beta)$, and hence is linear in $V$, and
\[
|c_\beta(V)| \leq c(\Phi) \, \norm{V}_{(M)} \, \norm{ f_\beta }_{(N)} .
\]
Thus there exists $C_\beta$ in $X_{(M)}'$ such that $c_\beta(V) = C_\beta(V)$.
\end{proof}

\section{Operators}\label{s:Operators}
We consider operators  $K$ on $L^2(\bR^n)$ defined by locally integrable kernels $k$, that is,
$$
Kf(x) = \int_{\bR^n} k(x,y) \, f(y) \,dy
\qquad\forall x \in \bR^n
$$
for all $f$ in $L^2(\bR^n)$.   Such operators may be given different orderings.  We write $|K_{1}| \leq |K_{2}|$ if
\begin{equation}\label{e:spectral-order}
 \norm{  K_{1} f }_{2} \leq \norm{  K_{2} f }_{2}
\qquad\forall f \in L^2(\bR^n)
\end{equation}
(this means that $K_{1}^*K_{1} \leq K_{2}^*K_{2}$ in the usual ordering of self-adjoint operators) and $|k_{1}| \leq |k_{2}|$ if
\begin{equation*}
|k_{1} (x,y)| \leq |k_{2} (x,y)|
\qquad\forall x,y \in \bR^n.
\end{equation*}
We may omit some of the absolute value signs if one of the operators is positive, or if one of the kernels is positive.
Note that if $K_{2}$ has an (unbounded) inverse $K_{2}^{-1}$ with a dense domain, then the inequality \eqref{e:spectral-order} amounts to saying that $K_{1} K_{2}^{-1}$ extends to a bounded operator of norm at most one.

\subsection{The heat operator.}
The heat semigroup is an important family of operators on $L^2(\bR^n)$.  For positive $t$, we define the heat kernel $p_t$ (often described as ``a Gaussian'') by
$$
p_t(x,y) := \frac{1}{t^{n/2}} \, \exp\Bigl(- \, \frac{\pi |x-y|^2}{t}\Bigr)
\qquad\forall x, y \in \bR^n.
$$
We then define the heat operator $P_t$ to be the operator corresponding to this kernel.
Then $\sF P_t f = \exp(-t\pi \abs{\,\cdot\,}^2) \sF f$.

We define the Laplacian $\Delta$ as a positive operator:
\[
\Delta := -\frac{1}{4\pi} \, \sum_{i=1}^n \frac{\partial^2}{\partial x_i^2},
\]
so $\sF(\Delta f) = \pi \abs{\,\cdot\,}^2 \hat f$, and $P_t$ may be written as $\exp(-t\Delta)$.  The operator $P_t$ has an unbounded inverse~$\exp(t\Delta)$ with dense domain.

\subsection{The main results.}
Implicitly, we use three quadratic forms $B_{1}$, $B_{2}$ and $B_{3}$ on $\bR^n \times \bR^n$:
\[
B_{1}(x,y) = 2 x\cdot y \qquad B_{2}(x,y) = |y|^2 \qquad B_{3}(x,y) = |x-y|^2
\]
for all $x$ and $y$ in $\bR^n$.  Recall that we define $E_{\rm i,1}(x,y) = \exp(2\pi i x\cdot y)$.  We also write $E_{\rm r,2}$ and $E_{\rm r,3}^t$ for the exponential functions $(x,y) \mapsto \exp(\pi\abs{y}^2)$ and $(x,y) \mapsto \exp(t \pi|x-y|^2)$ on $\bR^{2n}$.  The subscripts $\mathrm{r}$ and $\mathrm{i}$ stand for real and imaginary, and the subscripts $1$, $2$ and $3$ describe the quadratic form involved.  The superscript $t$ is used to indicate the real number $t$ in the exponent.

We are moving towards proving Theorem \ref{main-theorem-integrals}.  According to our discussions above, this is a theorem about operators $K$ and their locally integrable kernels $k$, such that $K \exp(t\Delta)$ is bounded on $L^2(\bR^n)$ and $E_{\rm r,3}^t k \in L^{\infty}(\bR^{2n})$.

The following version of Hardy's theorem for operators  boils down to Theorem \ref{Demangethm} if the operator $T$ from $\sS(\bR^n)$ to $\sS(\bR^n)'$ is given by a convolution,  but otherwise it is more general.

\begin{theorem}\label{thm:BDJ-op}
Suppose that $t$ is in $\bR^+$, and that $T$ is an operator from $\sS(\bR^n)$ to $\sS(\bR^n)'$ such that $T \exp(t\Delta)$ also maps $\sS(\bR^n)$ into $\sS(\bR^n)'$ and the kernel $\Phi$ in $\sS(\bR^{2n})'$ of $T$ (given by the Schwartz kernel theorem) satisfies $E_{\rm r,3}^t \Phi \in \sS(\bR^{2n})'$.
Then there exist a positive integer $N$ and distributions $\Mu_\alpha$ in $\sS(\bR^n)'$ when $\abs{\alpha} \leq N$ such that
\[
\Phi(x,y) = \sum_{\abs{\alpha}\leq N} \Mu_\alpha(x) \, (x-y)^\alpha  p_{t}(x,y)
\qquad\forall x,y \in \bR^n.
\]
\end{theorem}

\begin{proof}
By rescaling, we may suppose that $t=1$.  By our hypotheses,
$$\begin{aligned}
E_{\rm r,3}^1 \Phi  & \in \sS(\bR^{2n})' \\
E_{\rm r,2} \sF_{2} \Phi  & \in \sS(\bR^{2n})'.
\end{aligned}$$

Define the distribution $\Psi$ by 
\[
\Psi(x,y) = \Phi(x, x-y)
\qquad\forall x, y \in \bR^n.
\]
Then $E_{\rm r,2} \Psi(x,y) = E_{\rm r,3}^1 \Phi(x, x-y) \in \sS(\bR^{2n})'$, while
\[\begin{aligned}
E_{\rm r,2} \sF_{2}\Psi(x, \eta)
&= \exp(\pi|\eta|^2) \int_{\bR^n} \exp(-2\pi i \eta\cdot y) \, \Psi(x,y) \,dy \\
&= \exp(\pi|\eta|^2) \int_{\bR^n} \exp(-2\pi i \eta\cdot y) \, \Phi(x,x-y) \,dy \\
&= \exp(\pi|\eta|^2) \, E_{\rm i,1} \sF_{2} \Phi(x,-\eta) \,dy .
\end{aligned}\]

By the earlier result about multiplication by the function $E_{\rm i,1}$ and the hypothesis, $E_{\rm r,2} \sF_{2}\Psi \in \sS(\bR^{2n})'$.  The conclusion now follows from Theorem \ref{thm:BDJ-vv}, with $X$ taken to be $\sS(\bR^n)$.
\end{proof}

We are now able to prove Theorem \ref{main-theorem-integrals}, which we recall for the reader's convenience:
suppose that $t$ is in $\bR^+$, that $K$ is the operator on $L^2(\bR^n)$ defined by a locally integrable kernel $k$, and that
$$\begin{aligned}
|k|   &\leq  p_{t} \\
|K|  & \leq  P_{t}.
\end{aligned}$$
Then there exists a bounded measurable function $m$ on $\bR^n$ such that
\[
k(x,y) = m(x) \, p_{t}(x,y)
\qquad\forall x,y \in \bR^n.
\]

If $m$ is a bounded measurable function on $\bR^n$ and $k(x,y) = m(x)\,p_{t}(x,y)$, then $K$ is the heat operator $P_{t}$ followed by multiplication by $m$; the two inequalities of the theorem hold, but clearly $k$ may be more general than a Gaussian.

\begin{proof}[Proof of Theorem \ref{main-theorem-integrals}]
By rescaling if necessary, we may suppose that $t=1$.

For a nonnegative integer $j$, write $j'$ for $\max\{0, j-1\}$. For $f$ and $g$ in $\sS(\bR^n)$,
\begin{align*}
\Bigabs{ \bigl\langle f, K \Delta^j g \bigr\rangle }
&\leq \norm{ f }_{2}  \, \norm{  K \Delta^j g }_{2} \\
&\leq c\,\norm{ f }_{2} \, \norm{  P_{1} \Delta^j g }_{2}  \\
& =    c\,\norm{ f }_{2} \, \norm{  P_{1} \Delta^{j'} \Delta^{j-j'} g }_{2}  \\
& \leq    c\,\norm{ f }_{2} \,\sup \bigl\{  \exp( -\pi \abs{\xi}^2) \, (\pi |\xi|^{2})^{j'} : \xi \in \bR^n \bigr\} \,
\bignorm{\Delta^{j-j'} g}_{2} \\
& \leq    c\, \biggl(\frac{j'}{e}\biggr)^{j'} \norm{ f }_{(0)} \, \norm{ g }_{(2)}
\end{align*}
(when $j' = 0$, the right hand side is interpreted as $c\norm{ f }_{(0)}\norm{ g  }_{(2)}$).

Now
 \begin{align*}
 \abs{ \langle f, K \exp(\Delta) g \rangle }
&\leq \sum_{j\in\bN} \frac{1}{j!} \, \bigabs{ \left\langle f, K \Delta^j g\right\rangle } \\
&\leq   c\, \sum_{j \in\bN} \frac{1}{j!} \biggl(\frac{j'}{e}\biggr)^{j'} \norm{ f }_{(0)} \, \norm{ g }_{(2)}  \\
&\leq c\,\ \norm{ f }_{(0)} \, \norm{ g }_{(2)}
\qquad\forall f, g \in \sS(\bR^n) .
\end{align*}

Thus the kernel of the operator $K \exp(\Delta)$ is in $\sS(\bR^{2n})'$. Application of $\sF_{2}$ shows that $E_{\rm r,2} \sF_{2} k \in \sS(\bR^{2n})'$.

By Theorem \ref{thm:BDJ-op}, there exist $\Mu_\alpha$ in $\sS(\bR^n)'$ such that
\[
k(x,y) = \sum_{\abs{\alpha}\leq N} \Mu_\alpha(x) \, (x-y)^\alpha p_1(x,y)
\qquad\forall x,y \in \bR^n .
\]
But $|k| \leq p_{1}$, so $N=0$ and $\Mu_0 = m_0$ for some $m_0$ in $L^\infty(\bR^n)$.   \end{proof}

\end{document}